%% file: main.tex
\documentclass[preprint,review,11pt]{elsarticle}

\usepackage{natbib}
 \bibpunct[, ]{(}{)}{,}{a}{}{,}%
\usepackage{comment}
\usepackage{subfigure}
\usepackage{floatrow}
\usepackage{lscape}
\usepackage{newunicodechar}
\newunicodechar{✓}{\checkmark}
\usepackage{enumerate}
\usepackage{mathrsfs}
\usepackage[justification=centering]{caption}
\usepackage{algorithm,setspace}
\usepackage[noend]{algpseudocode}
\usepackage{epstopdf}
\usepackage{hyperref}
\usepackage{amsmath}
\usepackage{amssymb}
\usepackage{mathtools}
\usepackage{cases}

\makeatletter
\def\@author#1{\g@addto@macro\elsauthors{\normalsize%
    \def\baselinestretch{1}%
    \upshape\authorsep#1\unskip\textsuperscript{%
      \ifx\@fnmark\@empty\else\unskip\sep\@fnmark\let\sep=,\fi
      \ifx\@corref\@empty\else\unskip\sep\@corref\let\sep=,\fi
      }%
    \def\authorsep{\unskip,\space}%
    \global\let\@fnmark\@empty
    \global\let\@corref\@empty  
    \global\let\sep\@empty}%
    \@eadauthor={#1}
}
\makeatother

\journal{TBD}

\begin{document}

\begin{frontmatter}



\title{An Integrated Vaccination Site Selection and Dose Allocation Problem with Fairness Concerns}
\author{Linda Li\corref{cor1}\fnref{fn1}}
\cortext[cor1]{Corresponding Author; E-mail: jacksonli@missouristate.edu; Phone: (205)765-4530}
\author{Mohammad Firouz\fnref{fn2}}
\author{Daizy Ahmed\fnref{fn3}}
\author{Abdulaziz Ahmed\fnref{fn4}}

\address[fn1]{Department of Marketing, Missouri State University, Springfield, MO 65897.}
\address[fn2]{Department of Management, Information Systems, and Quantitative Methods, University of Alabama at Birmingham, Birmingham, AL 35294.}
\address[fn3]{Department of Mathematics and Statistics, Sri Mayapur International School, Sri Mayapur Dhama, Nadia District, West Bengal 741313.}
\address[fn4]{Department of Health Services Administration, University of Alabama at Birmingham, Birmingham, AL 35233.}

\begin{abstract}

\noindent In this paper we discuss \dots
\end{abstract}


\end{frontmatter}


\pagebreak

\input{1-Intro.tex}
\input{2-LitReview.tex}

\input{3-Model.tex}
\input{4-ModAnasis.tex}
\input{5-NumStudy.tex}
\input{6-Conclusions.tex}

\input{99-Appendices}

\newpage
\bibliographystyle{ormsv080}
\bibliography{main}

\end{document}

%% file: 1-Intro.tex
\section{Introduction}
\label{sec:Intro}
The ongoing Coronavirus Disease 2019 (COVID-19) pandemic has claimed the lives of over 4.7 Million people across the globe thus far and continues to do so under various emerging variants. Various vaccines have been dispatched into production and distribution in different countries to prevent the spread of the COVID-19, with the majority following a two-dose vaccination procedure. While 44.5\% of the world population has received at least one dose of a COVID-19 vaccine, this number amounts to only 2.2\% in the low-income countries, leaving majority of the population in high risk areas of the world in need of vaccination \citep{WHO}. 

Iran as one of the developing low-income countries in the Middle East is attempting on the distribution of vaccine to its population with the critical task of improving the immunity of its people towards COVID-19. Iranian Research Organization for Science and Technology (IROST) works directly under the Ministry of Science, Research, and Technology to develop a concrete plan to address this task.   

On the country level, the vaccine doses are first allocated to different states. Subsequently, each state allocates the doses of vaccine to major hospitals, where each hospital potentially serves several demand zones. To manage the distribution of vaccine doses effectively, the state mandates that each demand zone receive vaccine from their assigned hospital. Given the particular condition under which the vaccine can be stored and administered by the staff nurses, each hospital has a different capacity for vaccination. Additionally, each hospital located in a particular geographical location has different administration costs. Specifically, for hospital to serve as a vaccination site a hospital-specific fixed and variable cost is incurred. The total vaccination budget that the ministry receives is limited. 
While the decisions on the country level are taken by higher authority, IROST is interested in an efficient, equitable, and accessible dose allocation program in the state level. Specifically, the critical goals in the vaccination program is not only to have as much of the population vaccinated, but also to develop a plan that ensures fairness in distribution of the vaccine to the public as well as their accessibility in terms of travel distance to their assigned vaccination hospital.

\cite{beamon2008performance} define two pivotal performance measures to gauge humanitarian relief operations which also apply to the case of vaccine distribution: fairness and efficiency. In this paper, we address the problem of vaccination site selection and dose allocation by developing policies that are (i) efficient, (ii) equitable, (iii) accessible.

Given the scarcity of the available COVID-19 vaccines especially in a low-income country such as Iran, the efficiency in allocation of the vaccines is extremely important. Therefore, one of IROST's main goals is to achieve the highest level of vaccination among the public. Therefore, efficiency in our setting pertains to the percentage of the total supply of vaccine doses that are administered. 

Fairness in vaccination is not only important from a social justice point of view, but IROST's experience has shown that a fair distribution of vaccine proves more effective in public immunization by preventing highly-concentrated infected areas to form among the population. In this paper, we address fairness from two simultaneous points of view: equity and accessibility. Equity in our setting means that as far as possible, each demand zone should receive a fair-share of the total doses available. On the other hand, accessibility means that as far as possible, each demand zone should have equal travel distance to access their assigned vaccination site. 

\begin{figure}[H]
	\centering
	\includegraphics[width=0.3\textwidth]{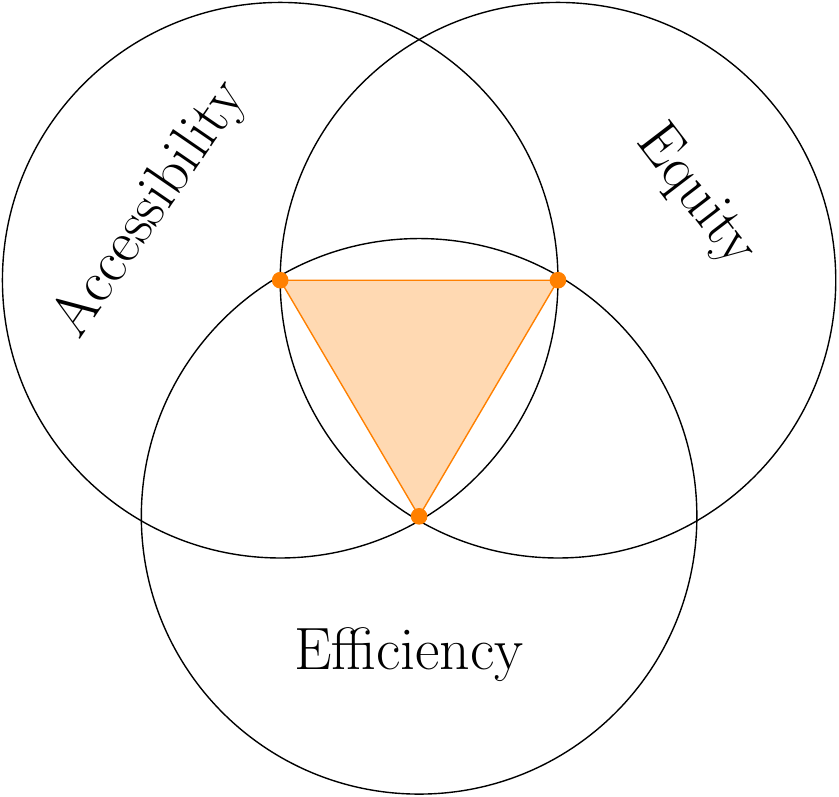}
	\caption{The dynamics of equitable, efficient, and accessible policies.}
	\label{fig:Equity_Efficiency_Accessibility}
\end{figure}

Efficiency, equity, and accessibility are conflicting objectives. Figure \ref{fig:Equity_Efficiency_Accessibility} shows the dynamics of the three performance measures in a vaccine distribution setting. The center of each circle represents the perfectness when each measure is considered by itself while the highlighted triangle represents the area in which the decision-maker in our setting may choose to perform in, depending on the flexibility of their policy regarding each one of the three performance measures as well as practical considerations in their network. 

%% file: 2-LitReview.tex
\section{Relevant Literature}
\label{sec:LitReview}

%% file: 3-Model.tex
\section{Model Formulation}
\label{sec:Model}
We consider the problem faced by a decision-maker that wants to locate $n$ vaccination sites within a region that contains $m$ demand zones, and simultaneously allocate the limited doses of vaccine to these sites in the most efficient and equitable manner. Specifically, the integrated decision seeks maximizing the percentage of demand fulfilled in each site while at the same time maintaining fairness in terms of both the fill-rates across the sites and the traveling distance across the zones. 

Table \ref{tab:param} lists the notations along with their definitions that we introduced so far and additional ones throughout the rest of the paper. 
\renewcommand{\arraystretch}{1.2}
\begin{table}[H]
\begin{tabular}{ll}
\hline
\textbf{Parameters} & \textbf{Definitions} \\
\hline
$n$ & number of candidate sites for vaccination. 		\\	
$m$ & number of demand zones.	\\	
$\mathcal{I}$ & set of candidate sites for vaccination, $\mathcal{I}=\{1,\dots,n\}$. 	\\	
$\mathcal{J}$ & set of demand zones, $\mathcal{J}=\{1,\dots,m\}$.	\\	
$S$ & the total supply of vaccine doses available.	\\
$B$ & the total allocation budget available.	\\
$C_i$ & total doses of vaccine that can be administered at site $i$. \\
$v_{i}$ & unit dose allocation cost in candidate site $i$. \\
$f_i$ & fixed cost of selecting candidate site $i$ to be a vaccination site. \\
$D_{j}$ & demand for doses of vaccine at zone $j$.	\\
$a_{j}$ & x-coordinate of the center of zone $j$. \\
$b_{j}$ & y-coordinate of the center of zone $j$. \\
$a_{i}$ & x-coordinate of the center of candidate site $i$. \\
$b_{i}$ & y-coordinate of the center of candidate site $i$. \\
$d_{ij}$ & demand weighted distance between zone $j$ and candidate site $i$. \\
$d_{j}$ & standardized distance from zone $j$ to its assigned vaccination site.\\
$\beta_{j}$ & fill-rate of zone $j$, $\beta_{j} =\frac{\sum\limits_{i \in \mathcal{I}}y_{ij}}{D_j}$.\\
$\beta$ & overall fill-rate, $\beta =\frac{\sum\limits_{i \in \mathcal{I}}\sum\limits_{j \in \mathcal{J}}y_{ij}}{\sum\limits_{j \in \mathcal{J}}D_j}$.\\
$\alpha$ & penalty for deviations from the maximum travel distance $\hat{d}$.	\\
$\theta$ & penalty for deviations from the maximum fill-rate $\hat{\beta}$.	\\
\hline
\textbf{Decision Variables} & \textbf{Definitions} \\
\hline
$x_{i}$ & 1 if candidate site $i$ is selected for vaccination, 0 otherwise.  \\
$y_{ij}$ & doses of vaccine provided to zone $j$ by vaccination site $i$. \\
$z_{ij}$ & 1 if zone $j$ assigned to vaccination site $i$, 0 otherwise.  \\
$\hat{\beta}$ & maximum fill-rate across the zones, $\hat{\beta}=\max\limits_{j\in \mathcal{J}} \beta_j$. \\
$\hat{d}$ & maximum travel distance across the zones, $\hat{d}=\max\limits_{j\in \mathcal{J}} d_{j}$. \\
\hline
\end{tabular}
\caption{Summary of the mathematical notation.}
\label{tab:param}
\end{table}

Our mathematical formulation is as follows.
 
\begin{align}
Max \hspace{3pt} \beta - \sum_{j\in\mathcal{J}} \frac{d_{j}}{m} -\sum_{j\in\mathcal{J}} \frac{\theta(\hat{\beta}-\beta_j)}{m} -\sum_{j\in\mathcal{J}}  \frac{\alpha (\hat{d}-d_{j})}{m} \label{eq:ObjFnc}
\end{align}\setcounter{equation}{0} 
\begin{subnumcases}{s.t.}
   \sum_{i\in\mathcal{I}}\sum_{j \in \mathcal{J}}y_{ij} \leq S & \label{eq:Const1}\\
   \sum_{i\in \mathcal{I}} f_ix_i+ \sum_{i\in\mathcal{I}}\sum_{j\in\mathcal{J}}v_{i} y_{ij} \leq B & \label{eq:Const2} \\
   \sum_{j \in \mathcal{J}}y_{ij} \leq C_i & $\forall i\in \mathcal{I}$ \label{eq:Const1}\\
   y_{ij} \leq Sz_{ij} & $\forall i\in \mathcal{I}\, ,\, \forall j\in \mathcal{J}$\label{eq:Const1}\\
   z_{ij} \leq x_i &$\forall i\in \mathcal{I} \, , \, \forall j\in \mathcal{J}$  \label{eq:Const1}\\
   \sum_{i\in\mathcal{I}}z_{ij} = 1 & $\forall j\in \mathcal{J}$ \label{eq:Const3}\\
   d_{ij}=D_j\sqrt{(a_j - a_i)^2 + (b_j-b_i)^2} & $\forall i\in \mathcal{I}\, ,\, \forall j\in \mathcal{J}$ \label{eq:Const4}\\
   \bar{d}_{ij}=\frac{d_{ij}-\min\limits_{i\in\mathcal{I},j\in\mathcal{J}}d_{ij}}{\max\limits_{i\in\mathcal{I},j\in\mathcal{J}}d_{ij}-\min\limits_{i\in\mathcal{I},j\in\mathcal{J}}d_{ij}} & $\forall i\in \mathcal{I}\, ,\, \forall j\in \mathcal{J}$ \label{eq:Const4}\\
   d_j =\sum\limits_{i\in\mathcal{I}}z_{ij}\bar{d}_{ij} & $\forall j\in \mathcal{J}$ \label{eq:Const4}\\
   d_j \leq \hat{d} & $\forall j\in \mathcal{J}$ \label{eq:Const4}\\
   \beta_{j} =\frac{\sum\limits_{i \in \mathcal{I}}y_{ij}}{D_j} & $\forall j\in \mathcal{J}$ \label{eq:Const5}\\
   \beta_j \leq \hat{\beta} & $\forall j\in \mathcal{J}$ \label{eq:Const5}\\
   y_{ij}\in \mathbb{N}_0; \; x_i, z_{ij} \in \{0,1\}; 0\leq \hat{\beta},\hat{d}\leq 1& $\forall i\in \mathcal{I}, \, \forall j\in \mathcal{J},$ \label{eq:Const6}
\end{subnumcases}
where $\mathbb{N}_0$ is the set of all non-negative integers and $\beta_{j} =\frac{\sum\limits_{i \in \mathcal{I}}y_{ij}}{D_j}$ is defined as the fill-rate at zone $j$ with $\hat{\beta}$ being the maximum among $\beta_j$'s. Additionally, $d_j =\sum\limits_{i\in\mathcal{I}}z_{ij}\bar{d}_{ij}$ is the standardized traveling distance from zone $j$ to its assigned vaccination site, where $\bar{d}_{ij}=\frac{d_{ij}-\min\limits_{i\in\mathcal{I},j\in\mathcal{J}}d_{ij}}{\max\limits_{i\in\mathcal{I},j\in\mathcal{J}}d_{ij}-\min\limits_{i\in\mathcal{I},j\in\mathcal{J}}d_{ij}}$ and $d_{ij}=D_j\sqrt{(a_j - a_i)^2 + (b_j-b_i)^2}$. 

In our objective function represented in Equation \eqref{eq:ObjFnc}, the first term is maximizing the fill-rate at each site and thereby seeks efficiency, while the second and third terms penalize inequity in terms of fill-rate and travelling distance and thereby seek equity. In detail, the second term penalizes the deviations from the maximum fill-rate across the sites by parameter $\theta$, while the third term penalizes the deviations from the maximum traveling distance across the zones by parameter $\alpha$. Higher values for $\theta$ and $\alpha$ set by the decision-maker, reflect the aptitude towards stricter requirements for maintaining equity in the policy. 

Constraint \eqref{eq:Const1} sets the limitation on the total number of doses of vaccine that are available for allocation. On the left hand-side of constraint \eqref{eq:Const2}, $f$ represents the fixed cost associated with establishing a site. On the other hand, there is a variable cost associated with each dose of vaccine that is to be administered at site $i$ which is captured by $v_i$. Constraint \eqref{eq:Const2}, therefore, sets the budget limitation on the total fixed and variable costs across the network. Constraint \eqref{eq:Const3} states that each demand zone must be exactly assigned to one vaccination site. Constraints \eqref{eq:Const4} and \eqref{eq:Const5} correspond to the definition of the maximum traveling distance and maximum fill-rate, respectively. Finally, Constraint \eqref{eq:Const6} sets the bounds on the decision variables. 

%% file: 4-ModAnasis.tex
\section{Model Analysis}
\label{sec:ModAnalysis}
The original optimization problem introduced in Model \eqref{eq:ObjFnc} looks at all the combinations of site-zone assignments by the definition of the decision variable $\bf z$. However, in reality, the demand zones form regional clusters and the assignment of sites-clusters instead of individual sites-zones can significantly reduce the computational burden. This is the basic idea behind Section \ref{subsec:alpha-fair-k-means} discussed below.

\subsection{$\alpha$-Fair Weighted k-Means Clustering}
\label{subsec:alpha-fair-k-means}
In this section, we introduce a fairness-adjusted version of the weighted k-means clustering technique to replace the site-zone assignment with site-cluster assignment. We note that the $\bf z$ vector in Model \eqref{eq:ObjFnc} has a solution space of $(2n)^m$, while using the cluster-based assignment, our algorithm has the solution space of $(2n)^{|\mathcal{K}|}$ for $\bf z$, where $|\mathcal{K}|$ is cardinality of the set of clusters $\mathcal{K}$ and we have $|\mathcal{K}|<<m$. 

\begin{figure}[H]
	\begin{center}
		\subfigure[Site-zone assignment.]{\includegraphics[width=0.47\textwidth]{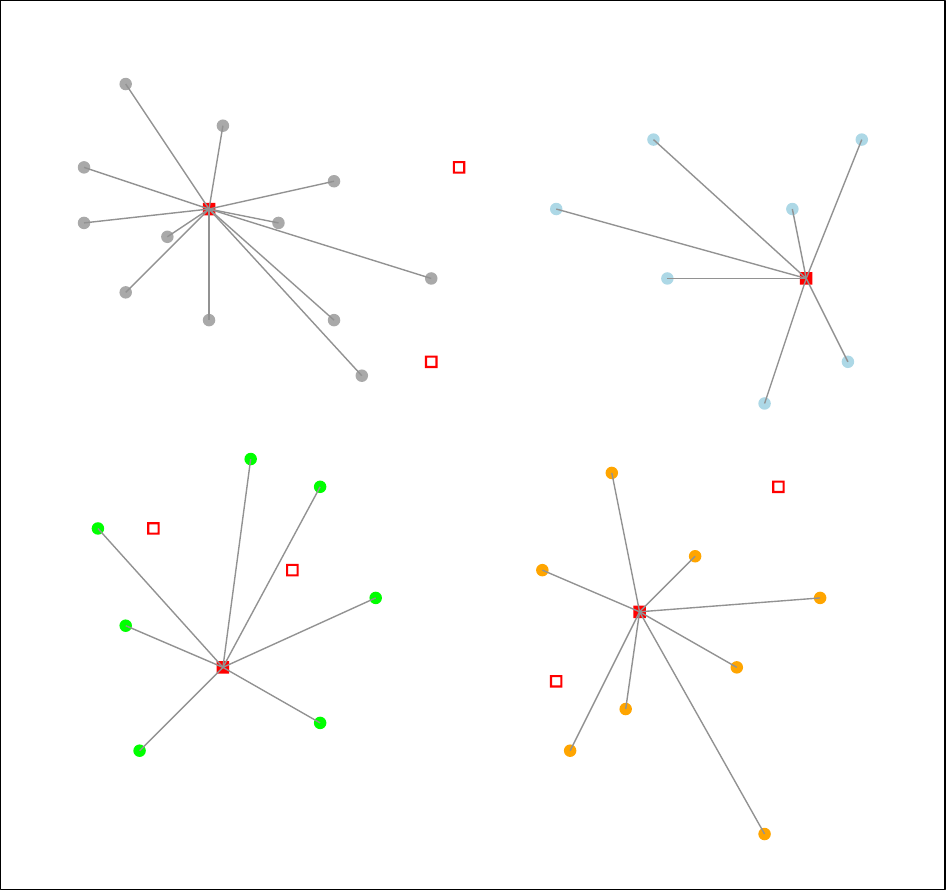} \label{fig:Site_Zone_Assignment_Example}}
		\subfigure[Site-cluster assignment.]{\includegraphics[width=0.47\textwidth]{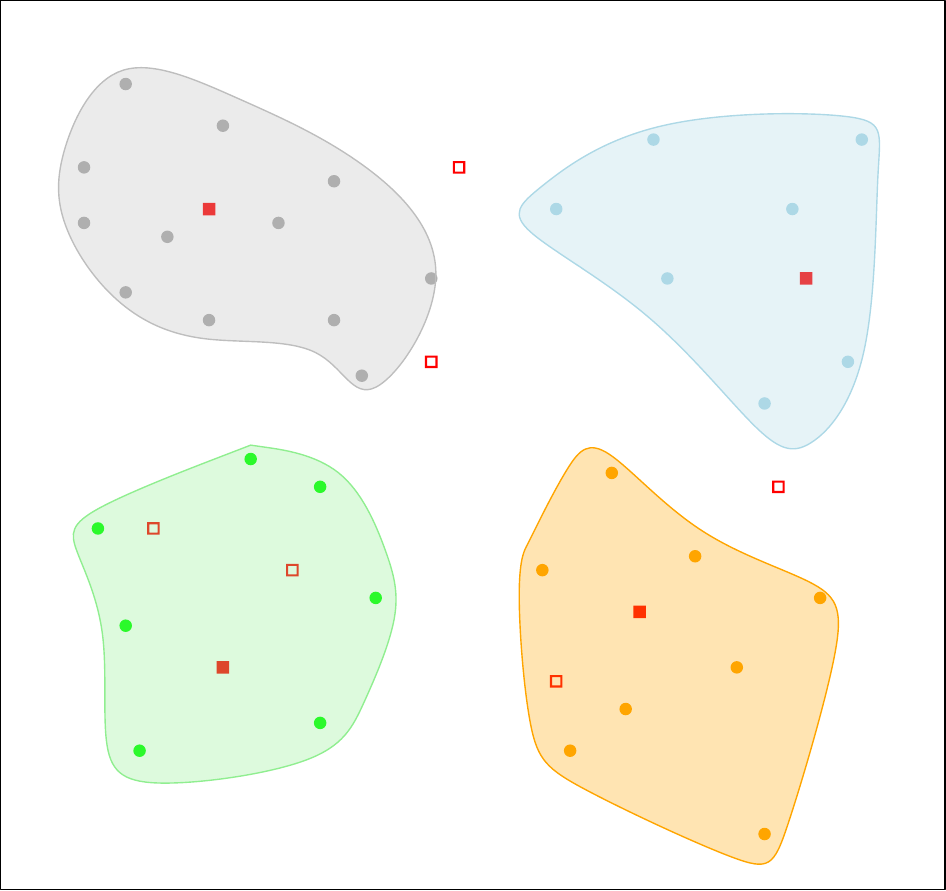} \label{fig:Site_Cluster_Assignment_Example}}
		\caption{problem size reduction through replacing site-zone with site-cluster assignment.}
		\label{fig:Site_Assignment}
	\end{center}
\end{figure}

Figure \ref{fig:Site_Assignment} demonstrates the idea with 34 demand zones and 10 candidate sites. The red squares represent the candidate sites while the various colored discs represent the demand zones in both Figures \ref{fig:Site_Zone_Assignment_Example} and \ref{fig:Site_Cluster_Assignment_Example}. In both Figures \ref{fig:Site_Zone_Assignment_Example} and \ref{fig:Site_Cluster_Assignment_Example}, the filled red squares are the selected sites to serve the demand zones for vaccination. In Figure \ref{fig:Site_Zone_Assignment_Example}, the solver has found the best values of $\bf z$ among $(2\times10)^{34}$ combinations of the solution space, while our algorithm has $(2\times10)^4$ combinations for $\bf z$ to deal with. We note that the selected sites are identical with both methods in this example. 

In Figure \ref{fig:Site_Cluster_Assignment_Example}, it is evident that the selected sites are those that fall closest to the centroid of the cluster. Exploiting this fact, next we specify how the size of the solution space can reduce further from $(2n)^{|\mathcal{K}|}$, by selecting a subset of promising candidate sites to consider for each cluster instead of all the sites.

\begin{enumerate}
    \item use weighted clustering (on both demand zones and candidate sites, with the weight of each zone being its demand for the vaccine and weight of each candidate site equal to one) to determine the global optimal number of clusters (with the highest Silhouette Coefficient). Call this set of clusters $\mathcal{K}$.
    \item if $n\geq |\mathcal{K}|$ continue, else change the number of clusters to $n$.
    \item set $e_{jk} =1$ if demand zone $j$ is assigned to cluster $k$
    \item make a list for each cluster in which candidate sites are sorted based on the following score $\frac{C_i}{d_{ik}f_iv_i}$, where $d_{ik}$ is the distance of site $i$ from the centroid of that cluster calculated as $\sqrt{(a_i-a_k)^2 +(b_i-b_k)^2 }$. Also, all the parameters $C_i$, $d_{ik}$, $f_i$, and $v_i$ are standardized to be between zero and one. 
    \item from each list in the previous step choose the first unique candidate site (start with the cluster that has the highest first score)
    \item within each cluster, capacity ($C_k$), fixed cost ($f_k$), and variable cost ($v_k$) become the corresponding values of the site assigned to that cluster.
    \item aggregate the demand in each cluster ($D_k$)
    \item allocate the total vaccine supply $S$ among the clusters by solving the following optimization problem:
        \begin{align}
            Max \hspace{3pt} \sum_{j\in\mathcal{J}}\sum_{k\in\mathcal{K}} \beta_{jk} -\sum_{j\in\mathcal{J}}\sum_{k\in\mathcal{K}} \theta(\beta-\beta_{jk})  \label{eq:ObjFnc_K}
            \end{align}\setcounter{equation}{0} 
            \begin{subnumcases}{s.t.}
            \sum_{j\in\mathcal{J}}\sum_{k\in\mathcal{K}}y_{jk} \leq S & \label{eq:Const1}\\
            \sum_{k\in\mathcal{K}}f_k+ \sum_{j\in\mathcal{J}}\sum_{k\in\mathcal{K}}v_{k} y_{jk} \leq B & \label{eq:Const2} \\
            \sum_{j\in\mathcal{J}}y_{jk} \leq C_k & $\forall k\in \mathcal{K}$ \label{eq:Const1}\\
            y_{jk} \leq e_{jk}S & $\forall j\in \mathcal{J}\, , \,\forall k\in \mathcal{K}$ \label{eq:Const5}\\
            \beta_{jk} = \frac{y_{jk}}{D_{jk}} & $\forall j\in \mathcal{J}\, , \,\forall k\in \mathcal{K}$ \label{eq:Const5}\\
            \beta_{jk} \leq \beta & $\forall j\in \mathcal{J}\, , \,\forall k\in \mathcal{K}$ \label{eq:Const5}\\
            y_{jk}\in \mathbb{N}_0; \; 0\leq \beta\leq 1& $\forall j\in \mathcal{J}\, , \,\forall k\in \mathcal{K}$ \label{eq:Const6}
        \end{subnumcases}
\end{enumerate}

%% file: 5-NumStudy.tex
\section{Numerical Study}
\label{sec:NumStudy}

%% file: 6-Conclusions.tex
\section{Conclusions}
\label{sec:Conclusions}

%% file: 99-Appendices.tex
\newpage
\begin{appendix}

\end{appendix}